\documentstyle[noteb,11pt,envcountreset]{lncse}

\font\elevenbb=msbm10 at 10.95pt

\def\C{\hbox{\elevenbb C}}
\def\N{\hbox{\elevenbb N}}
\def\M{\hbox{\elevenbb M}}

\def\Gr{Gr\"obner }
\def \bg #1 {\begin{tabular}{{#1}}}
\def \nd {\end{tabular}}
\def\Mma{{\it Mathematica\/}}
\def\exn #1#2 {\MathBegin{MathArray}[c]{l}
        #1 \\
        #2
  \MathEnd{MathArray}}
\newcommand \mwhile {{\bf while}\hspace{0.3cm}}

\newcommand \mdo {{\bf do}\hspace{0.3cm}}

\newcommand \mchoose {{\bf choose}\hspace{0.3cm}}

\newcommand \mbegin {{\bf begin}}
\newcommand \mend {{\bf end}}
\newcommand \bb {\hspace{0.3cm}}
\newcommand \h {\hspace{0.5cm}}
\newcommand \hh {\hspace{1.0cm}}
\newcommand \hhh {\hspace{1.5cm}}

\newcounter{cc}
\newcommand \hln {\hfill \addtocounter{cc}{1} \arabic{cc} \vskip 0.0cm \noindent }
\begin{document}
\title{Construction of Involutive Monomial Sets for Different
  Involutive Divisions}
\author{
  Vladimir P. Gerdt\inst{1} \and Vladimir V.Kornyak\inst{1}\and
  Matthias Berth\inst{2} \and G\"{u}nter Czichowski\inst{2}}
\institute{
  Laboratory of Computing Techniques and Automation\\
  Joint Institute for Nuclear Research\\
  141980 Dubna, Russia \\
  gerdt@jinr.ru \and
  Department of Mathematics and Informatics \\
  University of Greifswald \\
  D17487 Greifswald, Germany \\
  berth@rz.uni-greifswald.de \\
  czicho@rz.uni-greifswald.de
}
\maketitle

\begin{abstract}
  We consider computational and implementation issues for the
  completion of monomial sets to involution using different involutive
  divisions.  Every of these divisions produces its own completion
  procedure. For the polynomial case it yields an involutive basis
  which is a special form of a Gr\"{o}bner basis, generally redundant.
  We also compare our {\Mma} implementation of Janet division to an
  implementation in C.
\end{abstract}

\section{Introduction and Basic Definitions}
In our previous paper~\cite{GBC} we described our first results on
implementing in {\Mma} 3.0~\cite{Wolfram96} different involutive
divisions introduced in~\cite{GB1,GB2,G}; the completion of monomial
sets to involution for those divisions and application to constructing
Hilbert functions and Hilbert polynomials for monomial ideals.

In the present paper we pay more attention to efficient computation
and propose some algorithmic improvements. Besides, we shortly
describe an implementation of Janet division in C and compare the
running times for both implementations. Though in this paper we
consider involutivity of monomial ideals, all the underlying
operations with involutive divisions and monomials enter in more
general completion procedures for polynomial \cite{GB1,GB2} and
differential systems~\cite{Gerdt95}.

Let $\N$ be a set of non-negative integers, and
$\M=\{x_1^{d_1}\cdots x_n^{d_n}\ |\ d_i\in \N\}$
be a set of monomials in the polynomial ring $K[x_1,\ldots,x_n]$
over a field $K$ of characteristic zero.
By $deg(u)$ and $deg_i(u)$ we denote the total degree of $u\in \M$ and the
degree of variable $x_i$ in $u$, respectively. For the least common multiple
of two monomials $u,v\in \M$ we shall use the conventional notation
$lcm(u,v)$.
An admissible monomial ordering is denoted by $\succ$, and throughout this
paper we shall assume that it is compatible with
\begin{equation}
x_1\succ x_2\succ\cdots\succ x_n\,. \label{var_order}
\end{equation}

\begin{definition}{\cite{GB1}
An {\em involutive division} $L$ on $\M$ is given, if for any finite
monomial set $U\subset \M$ and for any $u\in U$ there is given a
submonoid $L(u,U)$ of $\M$ satisfying the conditions:
\begin{tabbing}
~~(a)~~\=If $w\in L(u,U)$ and $v|w$, then $v\in L(u,U)$. \\
~~(b)  \>If $u,v\in U$ and $uL(u,U)\cap vL(v,U)\not=\emptyset$, \\
       \>then $u\in vL(v,U)$ or $v\in uL(u,U)$. \\
~~(c)  \> If $v\in U$ and $v\in uL(u,U)$, then $L(v,U)\subseteq L(u,U)$. \\
~~(d)  \> If $V\subseteq U$, then $L(u,U)\subseteq L(u,V)$ for all $u\in V$.
\end{tabbing}
Elements of $L(u,U)$ are called {\em multiplicative} for $u$.  If
$w\in uL(u,U)$ we shall write $u|_L w$ and call $u$ an {\em
($L-$)involutive divisor} of $w$. In such an event the
monomial $v=w/u$ is {\em multiplicative} for $u$ and the equality
$w=uv$ will be written as $w=u\times v$. If $u$ is a conventional
divisor of $w$ but not an involutive one we shall write, as usual,
$w=u\cdot v$. Then $v$ is said to be {\em nonmultiplicative} for $u$.
} \label{inv_div}
\end{definition}

For every monomial $u\in U$, Definition~\ref{inv_div} provides the separation
\begin{equation}
\{x_1,\ldots,x_n\}=M_L(u,U)\cup NM_L(u,U),\\
\label{part}
\end{equation}
$M_L(u,U)\cap NM_L(u,U)=\emptyset$, of the set of variables into two subsets: {\em multiplicative}
$M_L(u,U)\subset L(u,U)$ and {\em nonmultiplica\-ti\-ve}
$NM_L(u,U)\cap L(u,U)=\emptyset$.
Conversely, if for any finite set $U\subset \M$ and any $u\in U$ the
separation~(\ref{part}) is given such that the corresponding
submonoid $L(u,U)$ of monomials in variables in $M_L(u,U)$ satisfies
the conditions (b)-(d), then the partition generates an involutive
division.

\begin{definition}{\cite{GB1}
Given an involutive division $L$, a
 monomial set $U$
 is {\em involutive} with respect to $L$ or $L-$involutive if
$$
 (\forall u\in U)\ (\forall w\in \M)\
 (\exists v\in U)\ \ [\ uw\in vL(v,U)\ ]\,.
$$
}
 \label{inv_mset}
\end{definition}

In this paper as well as in~\cite{GBC} we shall consider the following
eight different involutive divisions studied in~\cite{GB1,GB2,G}:

\begin{example}{Thomas division~\cite{Thomas}.
Given a finite set $U\subset \M$, the variable $x_i$ is considered as
multiplicative
for $u\in U$ if $deg_i(u)=max\{deg_i(v)\ |\   v\in U\}$, and
nonmultiplicative, otherwise.
} \label{div_T}
\end{example} 
\begin{example}{Janet division~\cite{Janet}. Let the set $U\subset \M$
be finite. For each $1\leq i\leq n$ divide $U$ into groups
labeled by non-negative integers $d_1,\ldots,d_i$:
$$[d_1,\ldots,d_i]=\{\ u\ \in U\ |\ d_j=deg_j(u),\ 1\leq j\leq i\ \}.$$
A variable $x_i$ is multiplicative for $u\in U$ if
$i=1$ and $deg_1(u)=max\{deg_1(v)\ |\ v\in U\}$,
or if $i>1$, $u\in [d_1,\ldots,d_{i-1}]$ and
$deg_i(u)=max\{deg_i(v)\ |\ v\in [d_1,\ldots,d_{i-1}]\}$.
} \label{div_J}
\end{example}

\begin{example}{
 Pommaret division~\cite{Pommaret78}. For a monomial
$u=x_1^{d_1}\cdots x_k^{d_k}$ with $d_k>0$ the variables $x_j,j\geq k$ are
considered as multiplicative and the other variables as nonmultiplicative.
For $u=1$ all the variables are multiplicative.
} \label{div_P}
\end{example}

\begin{example} {Division I~\cite{GB2}.
Let $U$ be a finite monomial set. The
variable
 $x_i$  is nonmultiplicative for $u\in U$ if there is $v\in U$ such
 that
 $$ x_{i_1}^{d_1}\cdots x_{i_m}^{d_m}u=lcm(u,v),\quad 1\leq m\leq [n/2],
 \quad d_j>0\ \ (1\leq j\leq m)\,,$$
 and $x_i\in \{x_{i_1},\ldots,x_{i_m}\}$.
} \label{div_I}
\end{example}

\begin{example} {Division II~\cite{GB2}.
For monomial
 $u=x_1^{d_1}\cdots x_k^{d_n}$ the variable $x_i$ is
 multiplicative if $d_i=d_{max}(u)$ where
 $d_{max}(u)=max\{d_1,\ldots,d_n\}$.
} \label{div_II}
\end{example}

\begin{example} {Induced division~\cite{G}.
Given an admissible monomial
 ordering $\succ $ a variable $x_i$
 is nonmultiplicative for $u\in U$ if there is $v\in U$ such that
 $v\prec u$ and $deg_i(u)<deg_i(v)$.
} \label{ind_div}
\end{example}

\noindent
To distinguish these divisions we use the abbreviations $T,J,P,I,II,D$
In the implementation described below, three orderings are used to induce
division in Example~\ref{ind_div}: lexicographical, degree-lexicographical
and degree-reverse-lexicographical. For these three induced divisions we shall
use the subscripts ${L}$, ${DL}$, ${DRL}$, respectively.

Every of the above divisions generates its own procedure for completion of
a monomial set to involution by means of its enlargement with involutively
irreducible nonmultiplicative prolongations. Given a monomial basis and an
involutive division, the following algorithm
{\bf MinimalInvolutiveMonomialBasis}~\cite{GB1} produces the uniquely defined
minimal involutive basis of the ideal.

\vskip 0.3cm
\noindent
\h Algorithm {\bf MinimalInvolutiveMonomialBasis:}
\vskip 0.2cm
\noindent
\h {\bf Input:}  $U$, a finite monomial set
\vskip 0.0cm
\noindent
\h {\bf Output:} $\bar{U}$, the minimal involutive basis of $Id(U)$
\vskip 0.0cm
\noindent
\h \mbegin
\hln
\hh $\bar{U}:=Autoreduce(U)$
\hln
\hh \mchoose any admissible monomial ordering $\prec$
\hln
\hh \mwhile exist $u\in \bar{U}$ and $x\in NM_L(u,\bar{U})$ s.t.
\hln
\hhh $u\cdot x$ has no involutive divisors in $\bar{U}$\bb \mdo
\hln
\hhh \mchoose such $u,x$ with the lowest $u\cdot x$ w.r.t. $\prec$
\hln
\hhh $\bar{U}:=\bar{U}\cup \{u\cdot x\}$
\hln
\hh \mend
\hln
\h \mend
\hln
\vskip 0.2cm
\noindent
Here $Autoreduce(U)$ stands for the conventional (non-involutive)
autoreduction.

\section{Implementation Issues}

In this section we will describe some observations that allow to speed
up the steps of the algorithm {\bf MinimalInvolutiveMonomialBasis}
significantly. Some of them are applicable to different divisions,
others are concerned with the completion procedure in general. The
basic operations on monomial sets are the same for the computation of
involutive bases of polynomial~\cite{GB1,GB2} and differential
systems~\cite{Gerdt95}, so the improvements described here are
relevant for these computations.

Our package provides a framework for studying the effect of using
different divisions and optimizations. It is implemented using a
``generic programming'' approach which allows to start with a straight
forward implementation of the algorithm and introduce more efficient
procedures for special situations later.

The following statement returns the minimal involutive basis of a
monomial set \inlineTFinmath{U} with respect to Janet division and
with lexicographic selection ordering:

{\center \dispSFinmath{\Muserfunction{minimalInvolutiveMonomialBasis}[\Mvariable{Janet}][U,\Mvariable{lexOrder}]}}

To extend the package for a new involutive division (called, say,
\mvb{newDivision}), one would only have to write the specific version
of the function \mvb{separation} which computes the multiplicative and
nonmultiplicative variables of a monomial $u\in U$ w.r.t. the set $U$:

{ \dispSFinmath{\Muserfunction{separation}[\Mvariable{newDivision}][\Mvariable{u\_},\Mvariable{U\_}]:=\ldots}}

\noindent All the other steps in the algorithm would then be executed by
functions that are generically defined for any involutive division.

On the other hand, an optimized procedure for a specific situation can
be introduced later to override the generic version. The pattern
matching mechanism in \Mma~ dispatches to the specific version
wherever it is appropriate.

Monomials are represented as multiindices, i.e. the
monomial \inlineTFinmath{x_{1}^{{i_1}}\cdot \ldots \cdot
  x_{n}^{{i_n}}} is represented as the list of its exponents
\inlineTFinmath{\{{i_1},\ldots ,{i_n}\}}. Thus, the set
\inlineTFinmath{U=\{{u_1},\ldots ,{u_m}\}} can be considered as a
\inlineTFinmath{m\times n-} matrix of integers. For every monomial
\inlineTFinmath{u}, we use two additional lists of length
\inlineTFinmath{n}: a list giving the separation of the variables for
$u$, and a similar list containing notes about the prolongations that
have already been done.

We will now describe observations that can be used to make the basic
operations of the algorithm {\bf MinimalInvolutiveMonomialBasis}
faster. Functions like {\it lcm\/} will be applied also to
multiindices, with the obvious meaning. The set notation is used for
lists, assuming that the order of the elements is given somehow.
\inlineTFinmath{U=\{{u_1},\ldots ,{u_n}\}} is a list of monomials, and
\inlineTFinmath{u} is always an element of \inlineTFinmath{U}.

The first step is to compute the separation for each of the input
monomials. For globally defined divisions, this is done irrespective
of the other monomials in \inlineTFinmath{U}.  For Janet division
(Example \ref{div_J}), we made use of the following remark:

\begin{remark}
  {When the list \inlineTFinmath{U} is sorted lexicographically in
    decreasing order, the groups \inlineTFinmath{[{d_1},\ldots
      ,{d_i}]} mentioned in the definition are grouped together. These
    groups are sorted lexicographically with respect to their labels
    of any fixed length $i$. The sorted list starts with the group
    labeled \inlineTFinmath{[{d_{1\Mvariable{max}}}],
      {d_{1\Mvariable{max}}}=\Mvariable{max} {{\Mvariable{deg}}_1}u},
    the monomials in \inlineTFinmath{[{d_{1\Mvariable{max}}}]} have
    \inlineTFinmath{{x_1}} as a multiplicative variable. We can split
    the list into groups given by labels of length 1 and proceed
    recursively within each of them, next considering degrees in the
    second variable \inlineTFinmath{{x_2}}, and so on.
\label{sep_J}
}
\end{remark}
\vspace{0.1cm}

For a division \inlineTFinmath{{D_{\succ }}} (Example \ref{ind_div})
that is induced by some ordering \inlineTFinmath{\succ }, we can use
an auxiliary list:
\begin{remark}
  {Let the monomials be sorted in descending order:
    \inlineTFinmath{{u_1}\succ \ldots \succ {u_n}}. We call the
    elements of the list
    \inlineTFinmath{\Muserfunction{cm}(U):=\{{m_1},\ldots ,{m_n} |
      {m_i}=\Muserfunction{lcm}({u_i},\ldots ,{u_n}),i=n,\ldots ,1\}}
    {\em the cumulated multiples} of \inlineTFinmath{U}. By
    definition, variable \inlineTFinmath{{x_j}} is nonmultiplicative
    for \inlineTFinmath{{u_i}} if and only if it has a higher degree
    in \inlineTFinmath{{m_i}}:
    \inlineTFinmath{{{\Mvariable{deg}}_j}{u_i}<{{\Mvariable{deg}}_j}{m_i}}.
    Thus, all we have to do is compute the list
    \inlineTFinmath{\Muserfunction{cm}(U) } of cumulated multiples and
    then compare each \inlineTFinmath{u\in U} against its
    corresponding entry in \inlineTFinmath{\Muserfunction{cm}(U) }.
\label{sep_ind_div}
}
\end{remark}

For Division I, we are not aware of any property that would allow
us to accelerate the computation of separations in a manner similar to
Janet or Induced divisions.

The following observation can be used to speed up the process of
finding a minimal nonmultiplicative prolongation (line 6 of the
algorithm).  Let us denote the minimal (w.r.t. the chosen ordering
\inlineTFinmath{\succ }) nonmultiplicative prolongation by a given
variable \inlineTFinmath{x} with \inlineTFinmath{{P_{\succ }}(x)}.

\begin{remark}
  {Let \inlineTFinmath{U} be sorted w.r.t. the completion
    ordering: \inlineTFinmath{{u_1}\succ \ldots \succ {u_n}}. Let
    \inlineTFinmath{{u_i} \Mvariable{and} x } be fixed such that
    \inlineTFinmath{{u_i}\cdot x} is a minimal nonmultiplicative
    prolongation w.r.t. \(\succ \). Then \inlineTFinmath{{u_i}\cdot x}
    is an element of the set \inlineTFinmath{\{{P_{\succ
          }}({x_1}),\ldots ,{P_{\succ }}({x_n})\}} .

    This follows directly from the minimality of
    \inlineTFinmath{{u_i}\cdot x}. Furthermore, \inlineTFinmath{{u_i}}
    is the minimal monomial having \inlineTFinmath{x} as a
    nonmultiplicative variable, because \inlineTFinmath{v\cdot x\succ
      u\cdot x} implies \inlineTFinmath{v\succ u}.}
\label{mini}
\end{remark}

\noindent The remark obviously extends to the more general situation of the
algorithm, where some of the nonmultiplicative prolongations have
already been considered. We keep a list
\inlineTFinmath{P = \{{P_{\succ }}({x_1}),\ldots
  ,{P_{\succ}}({x_n})\}} of nonmultiplicative prolongations, one for
each variable $x_1,\dots,x_n$, sorted by the completion ordering. Let
$v=u_i\cdot x_j$ be the minimal prolongation.  It is removed from $P$
and checked for involutive divisors. If $v$ is involutively reducible,
we have to add another prolongation w.r.t.  {\em the same variable
  $x_j$} to $P$. Otherwise, we add $v$ to the monomial set and
recompute the separations and $P$.

The next step in the algorithm is to search for an involutive divisor
\inlineTFinmath{w} of a nonmultiplicative prolongation
\inlineTFinmath{v=u\cdot x}. In the polynomial case, the efficiency of
this search can be even more important, since we may want to
involutively reduce every term of a prolonged polynomial.  Recall that
for an involutively reduced set \inlineTFinmath{U}, there can be at
most one such \inlineTFinmath{w}. We present now some optimizations
that apply to increasingly specialized situations.

The following remark uses a special property of involutive divisions,
taking into account that \inlineTFinmath{v} is a nonmultiplicative
prolongation of an element of \inlineTFinmath{U}.

\begin{remark}
  {Let \inlineTFinmath{U} be an involutively autoreduced set of
    monomials and \inlineTFinmath{v=u\cdot x} a nonmultiplicative
    prolongation of some \inlineTFinmath{u\in U}. If a monomial
    \inlineTFinmath{w\in U} is an involutive divisor of
    \inlineTFinmath{v} then
    \inlineTFinmath{{{\Mvariable{deg}}_x}w={{\Mvariable{deg}}_x}v}.

    Since \inlineTFinmath{u\cdot x} should be involutively reducible
    by \inlineTFinmath{w}, we can write \inlineTFinmath{u\cdot
      x=w\times (u\cdot x/w)}. If \inlineTFinmath{w=v=u\cdot x}, we
    are done. If \inlineTFinmath{w\neq u\cdot x} and
    \inlineTFinmath{w|u}, then \inlineTFinmath{u=w\times (u/w)}, which
    contradicts our assumption that \inlineTFinmath{U} is involutively
    autoreduced.  }
\label{finddiv}
\end{remark}

\noindent One can gain even more by considering particular divisions.
Consider a Janet-autoreduced set $U$. Let us denote the {\em longest
  common prefix} of two monomials $u,v$ by $lcp(u,v)$, where
$lcp(u,v):=(u_1,\dots,u_k)$ with $(u_1,\dots,u_k)=(v_1,\dots,v_k)$,
and $k$ the maximal index for which $u_k$ and $v_k$ coincide. If $u_1
\neq v_1$, we define $lcp(u,v):=()$. More generally, we use $lcp(v,U)$
to denote the longest common prefix that $v$ shares with some monomial
from the set $U$.

\begin{remark}
  {Assume that we search for a Janet - involutive divisor $w$ of a
    monomial $v$. Then, $w$ is in the class $\C$ defined by the label
    $lcp(v,U)$.
    Let $lcp(v,U)=(v_1,\dots ,v_k)$. Every involutive divisor
    $w=(w_1,\dots , w_n)$ is also a conventional divisor, thus $w_i
    \leq v_i, i=1,\dots ,k$. We show by contradiction that $w_i = v_i$
    for $i=1,\dots ,k$. Let $s$ be the smallest integer $1\leq s \leq
    k$ such that $w_s<v_s$. Then, $x_s$ is nonmultiplicative for $w$
    because there exists a monomial in the class $(v_1,\dots
    ,v_{s-1})$ which has higher degree in $x_s$, and $w$ is not an
    involutive divisor of $v$.  }
\label{finddiv_J}
\end{remark}

\noindent
Note that this remark applies to arbitrary monomials $v$, not only
those resulting from a nonmultiplicative prolongation.

Consider a nonmultiplicative prolongation $v=u\cdot x$. For Pommaret
division, an involutive divisor $w$ is reverse lexicographically
greater than $u$. For a division that is induced by
\inlineTFinmath{\succ }, either \inlineTFinmath{u\cdot x=w} or
\inlineTFinmath{u\succ w} holds.

These properties together with Remark~3.12 in~\cite{GBC}
suggest that one should keep the monomials sorted with respect to
some suitable order, and use this order as completion order, too.

Finally, when we find no involutive divisor, we have to add the
prolongation to the set and adjust separations for all monomials
accordingly. Except for globally defined divisions, this step is
potentially very time consuming.

\begin{remark}
  {For all divisions discussed so far, the following holds for a
    monomial \inlineTFinmath{u\in U}:
    \inlineTFinmath{\Muserfunction{NM}(u,U\cup
      \{v\})=\Muserfunction{NM}(u,U)\cup
      \Muserfunction{NM}(u,\{u,v\})}.  }
\label{pairwise}
\end{remark}
\noindent A detailed discussion of this fact can be found in \cite{G}.
After adding a monomial \inlineTFinmath{v} to \inlineTFinmath{U}, this
remark allows us to compute only the ``pairwise'' separations for
every \inlineTFinmath{u\in U}.

Specific divisions give rise to more improvements.

\begin{remark}
  {Let $v$ be a monomial, and assume that \inlineTFinmath{v} has
    no involutive divisor in the Janet-autoreduced set
    \inlineTFinmath{U}. Then, the separation may only change for
    monomials in the class $lcp(v,U)=(v_1, \dots ,v_k)$. The
    separation of the variables $x_1, \dots ,x_k$ is left unchanged.
    Furthermore, the separation of the variables $x_1, \dots ,x_k$ for
    the new monomial $v$ can be copied from the separation of any of
    the monomials in the class $lcp(v,U)$.  }
\label{rec_J}
\end{remark}

\begin{remark}
  {Consider adding a nonmultiplicative prolongation $v=u \cdot
    x_j$ to an autoreduced set w.r.t. some induced division
    \inlineTFinmath{{{\Mvariable{D}}_{\succ }}}.

    Only the variable \inlineTFinmath{{x_j}} can change from
    multiplicative to nonmultiplicative, and it can do so only for
    monomials \inlineTFinmath{s\succ v} satisfying
    \inlineTFinmath{{{\Mvariable{deg}}_j}s={{\Mvariable{deg}}_j}v-1}.
    }
\label{rec_ind}
\end{remark}

\noindent
Not all of the improvements mentioned here were actually implemented
in the package. Our experience suggests that sometimes the practical
performance in \Mma~ differs from what one expects from looking at the
algorithm. This is due to the interpreted nature of \Mma~ and its
flexible evaluation mechanism. Operations which are performed in the
kernel are usually much faster than their equivalent expressed in a
user defined function, and it was often a matter of trial and error to
decide which variant of an operation one should use for a given
division.

In practice, the size of the resulting involutive basis is certainly the
dominating factor for the overall running time of the algorithm. It was
thus worthwhile to invest more programming work in improvements for
those divisions which yield relatively small involutive bases (see
below).

The improvements for Janet division resulted the biggest gain in speed
compared to the generic implementation. When the completion ordering
is lexicographic, all optimizations described above are applied.
For induced divisions $D_\succ$, we always use $\succ$ as completion
ordering and Remark \ref{sep_ind_div} to recompute the separations.
Only for Division I, the time for changing the separations dominates
the time for the other basic operations. Division I is also the only
division for which the property mentioned in Remark \ref{pairwise} is
used.
The optimizations for finding an involutive divisor described above
have a positive effect for all divisions.

We have applied the package to examples taken from various sources.
For each polynomial system, we computed the degree - reverse -
lexicographical Gr\"{o}bner basis and took the resulting set of
leading monomials as input to the algorithm {\bf
  MinimalInvolutiveMonomialBasis}. As we described in~\cite{GBC} the
output can then be used to compute the Hilbert function, the Hilbert
polynomial and the index of regularity of the corresponding polynomial
ideal.

\begin{example}{
    \cite{BS92} Consider a \inlineTFinmath{n\times n} matrix
    \inlineTFinmath{A={{\big({{\alpha }_{i j}}\big)}_{n,n}}} with
    unspecified entries. The condition \inlineTFinmath{{A^2}=0} leads
    to a system of \inlineTFinmath{{n^2}} polynomial equations in the
    variables $\alpha_{11},\ldots ,\alpha_{1n},$ $ \alpha_{21},\ldots
    , \alpha_{nn}$. We treated the leading monomials of the degree
    reverse lexicographic Gr\"{o}bner basis, where the variables are
    ordered according to $\alpha_{11}\succ\ldots \succ
    \alpha_{1n}\succ \alpha_{21}\succ \ldots \succ \alpha_{nn}$.  }
\label{ex:sqn}
\end{example}

\begin{example}{
    The system of ``\inlineTFinmath{n} -th cyclic roots'' is a well
    known example. For \inlineTFinmath{n=4}, it is given by:}
\begin{eqnarray*}
        &&{x_1}+{x_2}+{x_3}+{x_4}=0, \\
        &&{x_1}{x_2}+{x_2}{x_3}+{x_3} {x_4}+{x_4} {x_1}=0, \\
        &&{x_1}{x_2}{x_3}+{x_2}{x_3}{x_4}+{x_3} {x_4}{x_1}+{x_4} {x_1}{x_2}=0, \\
        &&{x_1}{x_2}{x_3}{x_4}-1=0.
\end{eqnarray*}
\label{ex:cycn}
\end{example}

\begin{example}{The Reimer system in 5 variables:}
\begin{eqnarray*}
  && 1-2x_1^2+2x_2^2+2x_3^2+2*x_4^2-2x_5^2=0, \\
  && 1-2x_1^3+2x_2^3+2x_3^3+2*x_4^3-2x_5^3=0, \\
  && 1-2x_1^4+2x_2^4+2x_3^4+2*x_4^4-2x_5^4=0, \\
  && 1-2x_1^5+2x_2^5+2x_3^5+2*x_4^5-2x_5^5=0, \\
  && 1-2x_1^6+2x_2^6+2x_3^6+2*x_4^6-2x_5^6=0.
\end{eqnarray*}
\label{ex:r5}
\end{example}

\begin{example}{The Katsura system in 7 variables: }
\begin{eqnarray*}
&& x_1^2-x_1+2x_2^2+2x_3^2+2x_4^2+2x_5^2+2x_6^2+2x_7^2, \\ &&
2x_2x_1+2x_2x_3+2x_3x_4+2x_4x_5+2x_5x_6+2x_6x_7-x_2, \\ &&
2x_3x_1+2x_2x_4+2x_3x_5+2x_4x_6+2x_5x_7+x_2^2-x_3, \\ &&
2x_4x_1+2x_2x_5+2x_3x_6+2x_4x_7+2x_2x_3-x_4, \\ &&
2x_5x_1+2x_2x_6+2x_3x_7+2x_2x_4+x_3^2-x_5, \\ &&
2x_6x_1+2x_2x_7+2x_2x_5+2x_3x_4-x_6, \\ &&
x_1+2x_2+2x_3+2x_4+2x_5+2x_6+2x_7-1.
\end{eqnarray*}
\label{ex:k6}
\end{example}

\noindent
The following table shows the results of applying the algorithm {\bf
  MinimalInvolutiveMonomialBasis} to our examples. In the first three
columns, the size of the input is given where \inlineTFinmath{m} is
the number of monomials, \inlineTFinmath{n} is the number of
variables, and \inlineTFinmath{d} is the maximum total degree of the
input monomials. The divisions are indicated by the abbreviations used
above. For each division, we give the length of the minimal involutive
monomial basis, the number of prolongations considered during
completion, the portion of reducible prolongations, and the
computation time. Thus, 100$\%$ reducible prolongations means that the
input is already an involutive basis. An empty entry in the column for
Pommaret division means that we did not compute a minimal Pommaret
basis because the ideal is not zero dimensional. For the other
divisions, it means that the timing is larger than $10000$ seconds at
our computer\footnote{a 200 MHz 586 running Linux}.


\dispTFinmath{\MathBegin{MathArray}[c]{|l|c|c|c|c|c|c|c|c|c|c|c|}
  \hline \hline
  Input&\multicolumn{3}{c|}{Size}& \multicolumn{8}{c|}{Division}\\
  \cline{2-12}
  &m&n&d&J&T&P&I&\Mvariable{II}&{D_L}&{D_{\Mvariable{DRL}}}&{D_{\Mvariable{DL}}}
  \\
  \hline \hline

        {{\mathrm Ex.~} \ref{ex:r5}}&38&5&8&\MathBegin{MathArray}[c]{c}
        {\bf 55} \\
        190 \\
        91\% \\
        3.7\ s
  \MathEnd{MathArray}&\MathBegin{MathArray}[c]{c}
        {\bf 4392} \\
        17406 \\
        75\% \\
        4484\ s
  \MathEnd{MathArray}&\MathBegin{MathArray}[c]{c}
        {\bf 55} \\
        190 \\
        91\% \\
        3.4\ s
  \MathEnd{MathArray}&-&\MathBegin{MathArray}[c]{c}
        {\bf 151} \\
        503 \\
        77\% \\
        11\ s
  \MathEnd{MathArray}&\MathBegin{MathArray}[c]{c}
        {\bf 242} \\
        798 \\
        74\% \\
        48\ s
  \MathEnd{MathArray}&\MathBegin{MathArray}[c]{c}
        {\bf 894} \\
        3994 \\
        79\% \\
        556\ s
  \MathEnd{MathArray}&\MathBegin{MathArray}[c]{c}
        {\bf 594} \\
        2639 \\
        79\% \\
        267\ s
  \MathEnd{MathArray} \\ \hline
        {{\mathrm Ex.~} \ref{ex:k6}}&41&7&7&\MathBegin{MathArray}[c]{c}
        {\bf 43} \\
        211 \\
        99\% \\
        3.5\ s
  \MathEnd{MathArray}&-&\MathBegin{MathArray}[c]{c}
        {\bf 43} \\
        211 \\
        99\% \\
        3.7\ s
  \MathEnd{MathArray}&-&\MathBegin{MathArray}[c]{c}
        {\bf 201} \\
        861 \\
        81\% \\
        20\ s
  \MathEnd{MathArray}&\MathBegin{MathArray}[c]{c}
        {\bf 201} \\
        892 \\
        82\% \\
        44\ s
  \MathEnd{MathArray}&\MathBegin{MathArray}[c]{c}
        {\bf 1337} \\
        7600 \\
        83\% \\
        1500\ s
  \MathEnd{MathArray}&\MathBegin{MathArray}[c]{c}
        {\bf 1346} \\
        7663 \\
        83\% \\
        1539\ s
  \MathEnd{MathArray} \\ \hline
        \Mvariable{cyc}\ 4&7&4&6&\MathBegin{MathArray}[c]{c}
        {\bf 7} \\
        14 \\
        100\% \\
        0.19\ s
  \MathEnd{MathArray}&\MathBegin{MathArray}[c]{c}
        {\bf 98} \\
        242 \\
        62\% \\
        5.4\ s
  \MathEnd{MathArray}&-&\MathBegin{MathArray}[c]{c}
        {\bf 98} \\
        242 \\
        62\% \\
        18\ s
  \MathEnd{MathArray}&\MathBegin{MathArray}[c]{c}
        {\bf 25} \\
        55 \\
        67\% \\
        0.87\ s
  \MathEnd{MathArray}&\MathBegin{MathArray}[c]{c}
        {\bf 41} \\
        92 \\
        63\% \\
        2.3\ s
  \MathEnd{MathArray}&\MathBegin{MathArray}[c]{c}
        {\bf 9} \\
        20 \\
        90\% \\
        0.33\ s
  \MathEnd{MathArray}&\MathBegin{MathArray}[c]{c}
        {\bf 7} \\
        14 \\
        100\% \\
        0.21\ s
  \MathEnd{MathArray} \\ \hline
        \Mvariable{cyc}\ 5&20&5&8&\MathBegin{MathArray}[c]{c}
        {\bf 23} \\
        76 \\
        96\% \\
        1.1\ s
  \MathEnd{MathArray}&\MathBegin{MathArray}[c]{c}
        {\bf 1010} \\
        3544 \\
        72\% \\
        266\ s
  \MathEnd{MathArray}&\MathBegin{MathArray}[c]{c}
        {\bf 23} \\
        76 \\
        96\% \\
        1.1\ s
  \MathEnd{MathArray}&\MathBegin{MathArray}[c]{c}
        {\bf 1010} \\
        3544 \\
        72\% \\
        1656  s
  \MathEnd{MathArray}&\MathBegin{MathArray}[c]{c}
        {\bf 93} \\
        297 \\
        75\% \\
        5.5\ s
  \MathEnd{MathArray}&\MathBegin{MathArray}[c]{c}
        {\bf 154} \\
        488 \\
        72\% \\
        21\ s
  \MathEnd{MathArray}&\MathBegin{MathArray}[c]{c}
        {\bf 135} \\
        548 \\
        79\% \\
        21\ s
  \MathEnd{MathArray}&\MathBegin{MathArray}[c]{c}
        {\bf 106} \\
        419 \\
        79\% \\
        14\ s
  \MathEnd{MathArray} \\ \hline
        \Mvariable{cyc}\ 6&45&6&9&\MathBegin{MathArray}[c]{c}
        {\bf 46} \\
        194 \\
        99\% \\
        3.2\ s
  \MathEnd{MathArray}&-&\MathBegin{MathArray}[c]{c}
        {\bf 46} \\
        194 \\
        99\% \\
        3.1\ s
  \MathEnd{MathArray}&-&\MathBegin{MathArray}[c]{c}
        {\bf 201} \\
        807 \\
        81\% \\
        19\ s
  \MathEnd{MathArray}&\MathBegin{MathArray}[c]{c}
        {\bf 385} \\
        1527 \\
        78\% \\
        123\ s
  \MathEnd{MathArray}&\MathBegin{MathArray}[c]{c}
        {\bf 841} \\
        4230 \\
        81\% \\
        586\ s
  \MathEnd{MathArray}&\MathBegin{MathArray}[c]{c}
        {\bf 972} \\
        4899 \\
        81\% \\
        754\ s
  \MathEnd{MathArray} \\ \hline
  \exn{{\mathrm Ex.~} \ref{ex:sqn}}{n=3}
&25&9&4&\MathBegin{MathArray}[c]{c}
        {\bf 56} \\
        239 \\
        87\% \\
        4.5\ s
  \MathEnd{MathArray}&-&-&-&\MathBegin{MathArray}[c]{c}
        {\bf 612} \\
        2972 \\
        80\% \\
        131\ s
  \MathEnd{MathArray}&\MathBegin{MathArray}[c]{c}
        {\bf 531} \\
        2920 \\
        83\% \\
        313\ s
  \MathEnd{MathArray}&\MathBegin{MathArray}[c]{c}
        {\bf 1711} \\
        9362 \\
        82\% \\
        2593\ s
  \MathEnd{MathArray}&\MathBegin{MathArray}[c]{c}
        {\bf 1479} \\
        8044 \\
        82\% \\
        2048\ s
  \MathEnd{MathArray} \\ \hline
  \exn{{\mathrm Ex.~} \ref{ex:sqn}}{n=4}
        &161&16&6&\MathBegin{MathArray}[c]{c}
        {\bf 1324} \\
        11836 \\
        90\% \\
        923\ s
  \MathEnd{MathArray}&-&-&-&-&-&-&- \\
\hline \hline
  \MathEnd{MathArray}}

\vspace{\baselineskip}
\noindent For some examples, bases for two different divisions may
coincide. For the fourth cyclic roots
(Example \ref{ex:cycn}), the bases for Thomas division and Division I,
as well as those for Janet division and the induced division
\inlineTFinmath{{D_{\Mvariable{DL}}}} coincide, respectively.

The computations with monomial sets should give at least some hint to
the performance of different divisions in the polynomial and
differential cases. From our experience, Janet division, generally,
and Induced divisions, sometimes, seem to be the most promising in
terms of prolongations that have to be considered. Pommaret division
-- even though it is not noetherian -- deserves further investigation,
because it is globally defined and rather ``compact'', too.

\section{Conclusion}

In addition to the above described implementation of different
involutive divisions in {\Mma} we implemented the completion algorithm
for Janet division (Example~\ref{div_J}) in C. In this case an input
monomial set is represented as an array of lexicographically ordered
multiindices and its completion is done with respect to the same
order. This choice of completion ordering was motivated by the
monotonicity of Janet division with respect to the lexicographical
order. The partial involutivity of an intermediate monomial set is
preserved in the course of completion and the time for recomputation
of the separations is minimized~\cite{G}.

The set of nonmultiplicative prolongations to be treated is also
represented as a lexicographically sorted array of multiindices, that
provides the simplest way to choose a minimal prolongation.  Every
time an irreducible nonmultiplicative prolongation occurs it is
inserted in the intermediate monomial set and its nonmultiplicative
prolongations are inserted in the prolongation set. The determination
of their position in the sorted arrays is performed using the binary
search algorithm. In so doing, the check of Janet reducibility of the
prolongation under consideration is done in the course of the position
determination. This is a rather straightforward procedure that makes
use of the partition into prefix-groups as defined in
Example~\ref{div_J}.

The C implementation was done in GNU C/C++ version 2.81 on a
100 MHz Pentium computer running Windows 95. The running times for
examples in the above table are less than 0.01 seconds, except
Example~\ref{ex:sqn} for $n=4$ which took about 5 seconds.

We plan to extend both {\Mma} and C codes to polynomial and then to
linear differential systems. Whereas the highly flexible and easily
extensible {\Mma} code allows one to experiment with different
involutive divisions, in the further development of the C code we are
going to restrict ourselves to Janet, Pommaret and may be Induced
divisions which are more preferable from the computational efficiency
point of view.

\section{Acknowledgement}
The contribution of two of the authors (V.P.G. and
V.V.K.) was partially supported by grant INTAS-96-0842 and grant from the Russian
Foundation for Basic Research No. 98-01-00101.

\end{document}